\documentclass[10pt, letterpaper]{amsart}
\usepackage{amsmath, amsthm, amssymb}
\usepackage{mathbbol}
\usepackage{txfonts}

\newtheorem{thm}{Theorem}[section]
\newcommand{\bt}{\begin{thm}}
\newcommand{\et}{\end{thm}}

\newtheorem{cor}[thm]{Corollary}
\newcommand{\bc}{\begin{cor}}
\newcommand{\ec}{\end{cor}}

\newtheorem{lem}[thm]{Lemma}
\newcommand{\bl}{\begin{lem}}
\newcommand{\el}{\end{lem}}

\newtheorem{prop}[thm]{Proposition}
\newcommand{\bp}{\begin{prop}}
\newcommand{\ep}{\end{prop}}

\newtheorem{defn}[thm]{Definition}
\newcommand{\bd}{\begin{defn}}      
\newcommand{\ed}{\end{defn}}

\newtheorem{rmrk}[thm]{Remark}
\newcommand{\br}{\begin{rmrk}}
\newcommand{\er}{\end{rmrk}}

\newcommand{\thmref}[1]{Theorem~\ref{#1}}
\newcommand{\secref}[1]{Section~\ref{#1}}
\newcommand{\lemref}[1]{Lemma~\ref{#1}}

\newcommand{\propref}[1]{Proposition~\ref{#1}}

\newcommand{\N}{\mathbb{N}}

\newcommand{\R}{\mathbb{R}}
\newcommand{\Z}{\mathbb{Z}}

\newcommand{\dist}{\operatorname{dist}}
\newcommand{\diam}{\operatorname{diam}}

\newcommand{\hm}{{\mathcal H}}
\newcommand{\lm}{{\mathcal L}}

\newcommand{\lip}{\operatorname{Lip}}
\newcommand{\mass}[2][]{{\mathbf M_{#1}}(#2)}

\newcommand{\form}{{\mathcal D}}        
\newcommand{\curr}{{\mathbf M}}         
\newcommand{\intcurr}{{\mathbf I}}      

\newcommand{\flatnorm}{{\mathcal F}}
\newcommand{\fillvol}{{\operatorname{Fillvol}}}
\newcommand{\fillrad}{{\operatorname{Fillrad}}}

\newcommand{\rstr}{\:\mbox{\rule{0.1ex}{1.2ex}\rule{1.1ex}{0.1ex}}\:}
\newcommand{\bdry}{\partial}

\newcommand{\spt}{\operatorname{spt}}
\newcommand{\ohne}{\backslash}

\begin{document}

\title[A compactness theorem \`a la Gromov]{Compactness for manifolds and integral currents with bounded diameter and volume}

\author{Stefan Wenger}

\address
  {Department of Mathematics\\
University of Illinois at Chicago\\
851 S. Morgan Street\\
Chicago, IL 60607--7045 }
\email{wenger@math.uic.edu}

\date{October 29, 2008}

\thanks{Partially supported by NSF grant DMS 0707009}


\begin{abstract}
By Gromov's compactness theorem for metric spaces, every uniformly compact sequence of metric spaces admits an isometric embedding into a
common compact metric space in which a subsequence converges with respect to the Hausdorff distance. Working in the class or oriented
$k$-dimensional Riemannian manifolds (with boundary) and, more generally, integral currents in metric spaces in the sense of Ambrosio-Kirchheim and replacing the Hausdorff distance with the filling volume or flat distance, we prove an analogous compactness theorem in which we
replace uniform compactness of the sequence with uniform bounds on volume and diameter.
\end{abstract}

\maketitle

\section{Introduction}
In \cite{Gromov-polynomial-growth}, Gromov proved the following important result: Every uniformly compact sequence $X_n$ of metric spaces has a subsequence which converges in the Gromov-Hausdorff sense to a compact metric space $X$.  Recall that a sequence  of compact metric spaces $X_n$ is said to be {\it uniformly compact} if $\sup_n\diam X_n<\infty$ and if there exists a function $N: (0,\infty)\to\N$ such that for every $\varepsilon>0$ and $n\in\N$,  $X_n$ can be covered by at most $N(\varepsilon)$ balls of radius $\varepsilon$. In his article, Gromov in fact constructs a compact metric space
$Z$ into which every $X_n$ isometrically embeds.\footnote{Given isometric embeddings $\varphi_n: X_n\hookrightarrow Z$ it then follows that a subsequence $\varphi_{n_j}(X_{n_j})$ converges in the Hausdorff sense to a compact subset of $Z$.} 

In the present article we drop the condition of uniform compactness and show in our main theorem that, assuming only a uniform upper bound on diameter and volume, an analog of Gromov's result still holds if the Hausdorff distance is replaced by the filling volume or flat distance between Riemannian manifolds and more generally between integral currents in metric spaces. Before stating our main result, \thmref{theorem:general-compactness-fillvol}, in full generality, we state a version for closed Riemannian manifolds. For this let $M$ and $M'$ be closed oriented $k$-dimensional Riemannian manifolds and let $Z$ be a metric space into which $M$ and $M'$ isometrically embedd as metric spaces. Then $M$ and $M'$ may be viewed as singular Lipschitz $k$-cycles in $Z$, denoted by $\Lbrack M\Rbrack$ and $\Lbrack M'\Rbrack$, and their filling volume distance in $Z$ is defined as the least volume of a singular Lipschitz $(k+1)$-chain in $Z$ with boundary $\Lbrack M\Rbrack-\Lbrack M'\Rbrack$. A special case of our main theorem can now be stated as follows:

\bt\label{theorem:intro-manifold-case}
 Let $M_n$ be a sequence of closed oriented $k$-dimensional Riemannian manifolds with a uniform upper bound on diameter and volume. Then there exists a metric space $Z$, a subsequence $M_{n_j}$ and isometric embeddings $\varphi_j: M_{n_j}\hookrightarrow Z$ such that $\varphi_j(M_{n_j})$ is a Cauchy sequence with respect to the filling volume distance in $Z$ and converges to a  `generalized' Lipschitz $k$-cycle $T$ in $Z$. If $\varphi_j(M_{n_j})$ converges in the Hausdorff sense to a closed subset $Y\subset Z$ then $T$ is supported in $Y$.
\et 

By `generalized' Lipschitz cycle we in fact mean {\it integral current} in $Z$, see below.
Without imposing uniform bounds on diameter and volume the theorem is wrong: Let $M$ be a closed oriented $k$-dimensional Riemannian manifold. If $M_n$ is constructed by joining two copies of $M$ by a very thin tube of length $l_n\to\infty$, then $M_n$ cannot form a Cauchy sequence in any $Z$. Likewise, if $M_n$ is obtained from joining $n$ copies of $M$ by tubes of length $1$ to a fixed copy of $M$, then $M_n$ cannot form a Cauchy sequence either. 
Note that $M_n$ in the theorem need not have a Gromov-Hausdorff convergent subsequence in general. An example is given by gluing a sequence of thinner and thinner hairs to a fixed $M$ as above. An analog of \thmref{theorem:intro-manifold-case} holds for manifolds with boundary if the filling volume distance is replaced by the flat distance, see the main result below.

A natural framework in which to formulate our main theorem is provided by the theory of integral currents in complete metric spaces, developed by 
Ambrosio-Kirchheim in \cite{Ambr-Kirch-curr}, which extends the classical Federer-Fleming theory \cite{Federer-Fleming}. We recall that, given $k\geq 0$ and a complete metric space $Z$, the space of integral $k$-currents in $Z$ is denoted by $\intcurr_k(Z)$. The mass of $T\in\intcurr_k(Z)$ is denoted by $\mass{T}$. If $k\geq 1$, the boundary $\bdry T$ of $T$ is an element of $\intcurr_{k-1}(Z)$. Finally, for a Lipschitz map $\varphi: Z\to Z'$, the push-forward $\varphi_\#T$ of $T$ under $\varphi$ is an element of $\intcurr_k(Z')$. We refer to \secref{Section:Prelim} for definitions and details.
The main result of this article can be stated as follows:

\bt\label{theorem:general-compactness-fillvol}
 Let $k\in\N$, $C,D>0$ and let $X_n$ be a sequence of complete metric spaces. Given $T_n\in\intcurr_k(X_n)$ with $\mass{T_n}+\mass{\bdry T_n}\leq C$ and $\diam(\spt T_n)\leq D$ for all $n\in\N$ then there exists a subsequence $T_{n_j}$, a complete metric space $Z$, an integral current $T\in\intcurr_k(Z)$ and isometric embeddings $\varphi_j: X_{n_j}\hookrightarrow Z$ such that $\varphi_{j\#}T_{n_j}$ converges to $T$ in the flat distance in $Z$. If $\bdry T_n=0$ for all $n\in\N$ then $\fillvol(T-\varphi_{j\#}T_{n_j})\to 0$
as $j\to\infty$.
\et

Here, given $T,T'\in\intcurr_k(Z)$ with $k\geq 0$, their flat distance  is given by $d_{\mathcal F}(T,T'):= \flatnorm(T-T')$ where
\begin{equation*}
  \flatnorm(S):=\inf\{\mass{U}+\mass{V}:S=U+\bdry V, U\in\intcurr_k(Z), V\in\intcurr_{k+1}(Z)\}.
\end{equation*}
Note that
\begin{equation*}
 \flatnorm(S)\leq\min\{\fillvol(S),\mass{S}\}
\end{equation*}
for all $S\in\intcurr_k(Z)$ and that convergence in the flat distance in particular implies weak convergence of $T_n$ to $T$, that is, pointwise convergence.
%
%
%

First applications of \thmref{theorem:general-compactness-fillvol}, aiming at the asymptotic geometry of metric spaces, are given in \cite{Wenger-asymptotic}.
Note that in the theorem, the support $\spt T_n$ of $T_n$ need not be compact for any $n\in\N$. Furthermore, the above examples show that the theorem fails without the assumptions on diameter and mass and that, viewed as metric spaces endowed with the metric of $X_n$, the sequence $\spt T_n$ need not have a Gromov-Hausdorff convergent subsequence. We furthermore remark that Hausdorff convergence does not imply weak convergence of a subsequence: If, for example, $Z$ is the unit ball of an infinite dimensional separable Banach space, one can construct a sequence $T_n\in\intcurr_k(Z)$ satisfying the bounds in the theorem for which $\spt T_n = Z$ for all $n\in\N$ and $T_n$ does not have a weakly convergent subsequence. However, in the setting of the theorem, if $Y_j:=\varphi_{j}(\spt T_{n_j})$ Hausdorff converges in $Z$ then we can easily show that
\begin{equation}\label{equation:flatlim-subset-Hlim}
 \spt T\subset {\lim}_HY_j,
\end{equation}
see \propref{proposition:ultra-Hausdorff-limits}.
Note that in this case, $T$ has in general much more `regularity' than $ {\lim}_HY_j$ in the sense that $T$ is concentrated on a countably $\hm^k$-rectifiable set whereas $ {\lim}_HY_j$ may be arbitrarily complicated, for example infinite dimensional. Recall that countably $\hm^k$-rectifiable sets can be covered, up to an $\hm^k$-negligible set, by countably many images of biLipschitz maps from $\R^k$. In joint work, Sormani and the author have recently studied in \cite{Sormani-Wenger-short} conditions which imply equality in \eqref{equation:flatlim-subset-Hlim}. In general, if $\spt T_n$ does not have a Gromov-Hausdorff convergent subsequence, we can still show that $\spt T$ isometrically embeds into the ultralimit $(\spt T_{n_j})_\omega$, for every  non-principal ultrafilter $\omega$ on $\N$, see \secref{Section:Prelim}.

We turn to our second theorem which gives `uniqueness' of flat limits in the following sense.

\bt\label{theorem:uniqueness-intro}
 Let $(X_n)$ and $(T_n)$ be sequences as in the theorem above. Suppose there exist complete metric spaces $Z$, $Z'$ and isometric embeddings $\varphi_n: \spt T_n\hookrightarrow Z$ and $\varphi'_n:\spt T_n\hookrightarrow Z'$ such that $\varphi_{n\#}T_n$ converges in the flat distance to some $T\in\intcurr_k(Z)$ and $\varphi'_{n\#}T_n$ to some $T'\in\intcurr_k(Z')$. 
 If $T\not=0$ then $T'\not=0$ and there exists an isometry $\Psi: \spt T \to \spt T'$ with $\Psi_\#T=T'$.
\et

Note that we do not make any compactness assumptions in \thmref{theorem:uniqueness-intro}. A different formulation of the theorem will be given in \secref{Section:Uniqueness}. 

Sormani and the author have recently defined and studied in \cite{Sormani-Wenger-long} a notion of intrinsic flat distance between oriented $k$-dimensional Riemannian manifolds and, more generally, between certain countably $\hm^k$-rectifiable metric spaces, which is inspired by Gromov's idea of Gromov-Hausdorff distance. Using this new framework developed in \cite{Sormani-Wenger-long}, one can arrive at an elegant reformulation of the Theorems \ref{theorem:intro-manifold-case}, \ref{theorem:general-compactness-fillvol}, and \ref{theorem:uniqueness-intro} above. See
 \cite{Sormani-Wenger-long} for details.

\bigskip

\thmref{theorem:general-compactness-fillvol} seems to be new except in the special case that the sequence $(\spt T_n)$ of supports is uniformly compact. In this case it follows from known results as follows: By Gromov's theorem above, there exist a compact metric space $Z$ and isometric embeddings $\varphi_n: \spt T_n\hookrightarrow Z$. By Ambrosio-Kirchheim's Theorems 5.2 and 8.5 in \cite{Ambr-Kirch-curr} a subsequence $\varphi_{n_j\#}T_{n_j}$ converges weakly to an integral current $T\in\intcurr_k(Z)$. After possibly replacing $Z$ by $l^\infty(Z)$, Theorem 1.4 in \cite{Wenger-flatconv} shows that
$\varphi_{n_j\#}T_{n_j}$ converges to $T$ in the flat distance. In the general case, the rough idea is to decompose $T_n$ into the sum of two currents $T'_n$ and $T''_n$ such that $\spt T'_n$ is uniformly compact and $T''_n$ is close to $0$ in the flat distance.

\section{Preliminaries}\label{Section:Prelim}

In this section we review some definitions and facts from metric geometry and geometric measure theory which are used throughout the paper.

\subsection{Gromov-Hausdorff distance and Kuratowski embedding}

Let $Z$ be a metric space. The Hausdorff distance $d_H(A,B)$ between two subsets $A$ and $B$ in $Z$ is the infimum of all $\varepsilon\geq 0$ such that $A$ is contained in the $\varepsilon$-neighborhood of $B$ and $B$ is contained in the $\varepsilon$-neighborhood of $A$. As mentioned in the introduction, the Gromov-Hausdorff distance between two metric spaces $(X,d)$ and $(X', d')$ is
\begin{equation*}
  d_{GH}(X, X'):= \inf d_{H}(\varphi(X), \varphi'(X')),
\end{equation*}
where the infimum is taken over all metric spaces $Z$ and isometric embeddings $\varphi: X\hookrightarrow Z$ and $\varphi': X'\hookrightarrow Z$. It is not difficult to show that it is enough to take the infimum over all metric spaces $Z$ of the form $(X\sqcup X', d)$, where $d$ is a metric which coincides with $d$ on $X$ and with $d'$ on $X'$.

Given a metric space $(X, d)$, we denote by $l^\infty(X)$ the Banach space of bounded functions on $X$ together with the supremum norm. For fixed $x_0\in X$, the map $\iota(x):= d(\cdot, x) - d(\cdot, x_0)$ defines an isometric embedding of $X$ into $l^\infty(X)$, called Kuratowski embedding. It is well-known that $l^\infty(X)$ is an injective space for any $X$ in the following sense: Given any other metric space $Y$, a subset $A\subset Y$, and a $\lambda$-Lipschitz map $\varphi: A \to l^\infty(X)$, there exists an extension $\bar{\varphi}: Y\to l^\infty(X)$ of $\varphi$ which is $\lambda$-Lipschitz.

\subsection{Currents in metric spaces}
The theory of integral currents in metric spaces was developed by Ambrosio and Kirchheim in \cite{Ambr-Kirch-curr} and provides a suitable 
notion of generalized surfaces in metric spaces. In the following we recall the definitions from \cite{Ambr-Kirch-curr} that are needed throughout this paper.

Let $(X,d)$ be a complete metric space and $k\geq 0$ and let $\form^k(X)$ be the set of $(k+1)$-tuples $(f,\pi_1,\dots,\pi_k)$ 
of Lipschitz functions on $X$ with $f$ bounded. The Lipschitz constant of a Lipschitz function $f$ on $X$ will
be denoted by $\lip(f)$.
\bd
A $k$-dimensional metric current  $T$ on $X$ is a multi-linear functional on $\form^k(X)$ satisfying the following
properties:
\begin{enumerate}
 \item If $\pi^j_i$ converges point-wise to $\pi_i$ as $j\to\infty$ and if $\sup_{i,j}\lip(\pi^j_i)<\infty$ then
       \begin{equation*}
         T(f,\pi^j_1,\dots,\pi^j_k) \longrightarrow T(f,\pi_1,\dots,\pi_k).
       \end{equation*}
 \item If $\{x\in X:f(x)\not=0\}$ is contained in the union $\bigcup_{i=1}^kB_i$ of Borel sets $B_i$ and if $\pi_i$ is constant 
       on $B_i$ then
       \begin{equation*}
         T(f,\pi_1,\dots,\pi_k)=0.
       \end{equation*}
 \item There exists a finite Borel measure $\mu$ on $X$ such that
       \begin{equation}\label{equation:mass-def}
        |T(f,\pi_1,\dots,\pi_k)|\leq \prod_{i=1}^k\lip(\pi_i)\int_X|f|d\mu
       \end{equation}
       for all $(f,\pi_1,\dots,\pi_k)\in\form^k(X)$.
\end{enumerate}
\ed
The space of $k$-dimensional metric currents on $X$ is denoted by $\curr_k(X)$ and the minimal Borel measure $\mu$
satisfying \eqref{equation:mass-def} is called mass of $T$ and written as $\|T\|$. We also call mass of $T$ the number $\|T\|(X)$ 
which we denote by $\mass{T}$.
The support of $T$ is, by definition, the closed set $\spt T$ of points $x\in X$ such that $\|T\|(B(x,r))>0$ for all $r>0$. 

Every function $\theta\in L^1(K,\R)$ with $K\subset\R^k$ Borel measurable induces an element of $\curr_k(\R^k)$ by
\begin{equation*}
 \Lbrack\theta\Rbrack(f,\pi_1,\dots,\pi_k):=\int_K\theta f\det\left(\frac{\partial\pi_i}{\partial x_j}\right)\,d\lm^k
\end{equation*}
for all $(f,\pi_1,\dots,\pi_k)\in\form^k(\R^k)$.
The restriction of $T\in\curr_k(X)$ to a Borel set $A\subset X$ is given by 
\begin{equation*}
  (T\rstr A)(f,\pi_1,\dots,\pi_k):= T(f\chi_A,\pi_1,\dots,\pi_k).
\end{equation*}
This expression is well-defined since $T$ can be extended to a functional on tuples for which the first argument lies in 
$L^\infty(X,\|T\|)$.

If $k\geq 1$ and $T\in\curr_k(X)$ then the boundary of $T$ is the functional
\begin{equation*}
 \bdry T(f,\pi_1,\dots,\pi_{k-1}):= T(1,f,\pi_1,\dots,\pi_{k-1}).
\end{equation*}
It is clear that $\bdry T$ satisfies conditions (i) and (ii) in the above definition. If $\bdry T$ also satisfies (iii) then $T$ is called a normal current.
By convention, elements of $\curr_0(X)$ are also called normal currents.

The push-forward of $T\in\curr_k(X)$ 
under a Lipschitz map $\varphi$ from $X$ to another complete metric space $Y$ is given by
\begin{equation*}
 \varphi_\# T(g,\tau_1,\dots,\tau_k):= T(g\circ\varphi, \tau_1\circ\varphi,\dots,\tau_k\circ\varphi)
\end{equation*}
for $(g,\tau_1,\dots,\tau_k)\in\form^k(Y)$. This defines a $k$-dimensional current on $Y$.
It follows directly from the definitions that $\bdry(\varphi_{\#}T) = \varphi_{\#}(\bdry T)$.

\sloppy
We will mainly be concerned with integral currents. We recall that an $\hm^k$-measurable set $A\subset X$
is said to be countably $\hm^k$-rectifiable if there exist countably many Lipschitz maps $\varphi_i :B_i\longrightarrow X$ from subsets
$B_i\subset \R^k$ such that
\begin{equation*}
\hm^k\left(A\ohne \bigcup \varphi_i(B_i)\right)=0.
\end{equation*}
\fussy

An element $T\in\curr_0(X)$ is called integer rectifiable if there exist finitely many points $x_1,\dots,x_n\in X$ and $\theta_1,\dots,\theta_n\in\Z\ohne\{0\}$ such
that
\begin{equation*}
 T(f)=\sum_{i=1}^n\theta_if(x_i)
\end{equation*}
for all bounded Lipschitz functions $f$.
 A current $T\in\curr_k(X)$ with $k\geq 1$ is said to be integer rectifiable if the following properties hold:
 \begin{enumerate}
  \item $\|T\|$ is concentrated on a countably $\hm^k$-rectifiable set and vanishes on $\hm^k$-negligible Borel sets.
  \item For any Lipschitz map $\varphi:X\to\R^k$ and any open set $U\subset X$ there exists $\theta\in L^1(\R^k,\Z)$ such that 
    $\varphi_\#(T\rstr U)=\Lbrack\theta\Rbrack$.
 \end{enumerate}
Integer rectifiable normal currents are called integral currents. The corresponding space is denoted by $\intcurr_k(X)$. In case $X=\R^N$ is Euclidean space, $\intcurr_k(X)$ agrees with the space of $k$-dimensional Federer-Fleming integral currents in $\R^N$. If $A\subset\R^k$ is a Borel set of finite measure and
finite perimeter then $\Lbrack\chi_A\Rbrack \in\intcurr_k(\R^k)$. Here, $\chi_A$ denotes the characteristic function. If $T\in\intcurr_k(X)$ and if $\varphi:X\to Y$ is a Lipschitz 
map into another complete metric space then $\varphi_{\#}T\in\intcurr_k(Y)$. Every oriented $k$-dimensional complete Riemannian manifold $M$ with finite volume and finite boundary volume gives rise to an element $\Lbrack M\Rbrack \in\intcurr_k(M)$. Moreover, every Lipschitz chain in a complete metric space $X$ can be viewed as an integral current in $X$. Recently, a variant of Ambrosio-Kirchheim's theory that does not rely on the finite mass axiom has been developed by Lang in \cite{Lang-currents}.

\subsection{Filling volume and embeddings into ultralimits}
Let $Z$ be a complete metric space and $k\geq 0$. The filling volume in $Z$ of an element $T\in\intcurr_k(Z)$ is defined as
\begin{equation*}
 \fillvol(T):= \inf\{\mass{S}: \text{ $S\in\intcurr_{k+1}(Z)$ with $\bdry S = T$}\},
\end{equation*}
where we agree on the convention that $\inf\emptyset = \infty$, that is if no $S\in\intcurr_{k+1}(Z)$ with $\bdry S=T$ exists, e.g.\ in the case that $\bdry T\not=0$.

We prove the following easy facts mentioned in the introduction.

\bp\label{proposition:ultra-Hausdorff-limits}
 Let $C>0$ and let $T_n\in\intcurr_k(Z)$ be a sequence satisfying
 $\diam(\spt T_n)\leq C$ and  $\mass{T_n}\leq C$ for all $n\in\N$, and, in case $k\geq 1$, also $\mass{\bdry T_n}\leq C$.
 If there exists $T\in\intcurr_k(Z)$ such that $T_n$ converges weakly to $T$ then we have:
 \begin{enumerate}
  \item For every non-principal ultrafilter $\omega$ on $\N$ there exists an isometric embedding of $\spt T$ into the $\omega$-ultralimit of $(\spt T_n)$.  
  \item If $\spt T_n$ converges in the Hausdorff sense to a closed subset $Y\subset Z$ then $\spt T\subset Y$.
 \end{enumerate}
\ep

For the definitions and properties of non-principal ultrafilters on $\N$ and of ultralimits of sequences of metric spaces see e.g.\ \cite{Bridson-Haefliger}.

\begin{proof}
 Given $z\in\spt T$ and $\varepsilon>0$ there exist, by \cite[Proposition 2.7]{Ambr-Kirch-curr}, Lipschitz functions $f, \pi_1,\dots, \pi_k$ on $Z$ with $\spt f \subset B(z,\varepsilon)$ and such that 
 $T(f,\pi_1,\dots,\pi_k)\not= 0$. By the same proposition and by the definition of weak convergence there exists for every $n$ sufficiently large $z_n\in\spt T_n$ such that $d(z_n, z)\leq\varepsilon$. This shows that for every $z\in\spt T$ there exists a sequence $z_n\in Z$ converging to $z$ such that $z_n\in\spt T_n$ for every $n\in\N$. Statement (i) readily follows from this. As for (ii), one easily checks that the map $\varphi(z):= [(z_n)]$ from $\spt T$ to the ultralimit of the sequence of metric spaces $(\spt T_n, d_Z)$ is isometric. Here, $d_Z$ denotes the metric of $Z$.
\end{proof}

\section{A decomposition theorem for integral currents}\label{section:decomposition-theorem}
The main result of this section, \thmref{theorem:suitable-decomposition}, is a crucial ingredient of the proof of  \thmref{theorem:general-compactness-fillvol}. In what follows, we work in somewhat greater generality than would be needed. This will allow us to prove new isoperimetric estimates in \cite{Wenger-asymptotic}.

\bd
 Let $k\geq 2$ and $\alpha>1$. A complete metric space $X$ is said to admit an isoperimetric inequality of rank $\alpha$ for $\intcurr_{k-1}(X)$ 
 if 
 there is a constant $D>0$ such that for every $T\in\intcurr_{k-1}(X)$ with $\bdry T=0$ there exists $S\in\intcurr_k(X)$ with $\bdry S=T$ and
 \begin{equation}\label{equation:isop-ineq-alpha}
 \mass{S}\leq DI_{k,\alpha}(\mass{T}),
 \end{equation}
 where $I_{k,\alpha}$ is the function given by
 \begin{equation*}
 I_{k,\alpha}(r):=\left\{\begin{array}{l@{\qquad}l}
  r^{\frac{k}{k-1}} & 0\leq r\leq 1\\
  r^{\frac{\alpha}{\alpha-1}} & 1 < r <\infty.
 \end{array}\right.
\end{equation*}
\ed
In \cite[6.32]{Gromov-metric-structures} the polynomial bound $r^{\frac{\alpha}{\alpha-1}}$ was termed an isoperimetric inequality 
{\it of rank greater than} $\alpha$.
Here we will use the shorter terminology {\it of rank} $\alpha$.
Isoperimetric inequalities of rank $k$ for $\intcurr_{k-1}(X)$ are called of {\it Euclidean type}. Every Banach space $X$ admits isoperimetric inequalities of Euclidean type for $\intcurr_{k-1}(X)$ for every $k\geq 2$, see \cite{Wenger-GAFA}.

Set $\Lambda:= \{(k,\alpha)\in\N\times(1,\infty): k\geq 2\}\cup\{(1,0)\}$, let $\gamma\in(0,\infty)$ and define auxiliary functions by
\begin{equation*}
 F_{1,0,\gamma}(r)=\gamma r\quad\text{and}\quad G_{1,0}(r)=r
\end{equation*}
and for $(k,\alpha)\in\Lambda\ohne\{(1,0)\}$
\begin{equation*}
 F_{k,\alpha,\gamma}(r):=\left\{\begin{array}{l@{\qquad}l}
  \gamma\cdot r^k & 0\leq r\leq 1\\
  \gamma\cdot r^\alpha & 1 < r <\infty
 \end{array}\right.
\end{equation*}
and
\begin{equation*}
 G_{k,\alpha}(r):=\left\{\begin{array}{l@{\qquad}l}
  r^{\frac{1}{k}} & 0\leq r\leq 1\\
  r^{\frac{1}{\alpha}} & 1 < r <\infty.
 \end{array}\right.
\end{equation*}

The following is the main result of this section.

\bt\label{theorem:suitable-decomposition}
 Let $X$ be a complete metric space, $(k,\alpha)\in\Lambda$, and suppose in case $k\geq 2$ that $X$ admits an isoperimetric inequality of rank $\alpha$ for 
 $\intcurr_{k-1}(X)$. Then for every $\lambda\in(0, 1)$ there exists a $\gamma\in(0,1)$ with the following property. 
 Abbreviate $F:= F_{k,\alpha,\gamma}$ and $G:= G_{k,\alpha}$ and let $\delta\in(0,1)$.
 For every $T\in\intcurr_k(X)$ there exist $R\in\intcurr_k(X)$ and $T_j\in\intcurr_k(X)$, $j\in\N$, such that
 \begin{equation*}
  T= R + \sum_{j=1}^\infty T_j
 \end{equation*}
 and for which the following properties hold:
 \begin{enumerate}
  \item $\bdry R=\bdry T$ and $\bdry T_j=0$ for all $j\in\N$;
  \item For all $x\in\spt R\ohne\spt\bdry T$ and all $0\leq r\leq \min\{5\delta G(\mass{R}), \dist(x,\spt\bdry T)\}$
   \begin{equation*}
    \|R\|(B(x,r))\geq \frac{1}{2}5^{-(k+\alpha)}F(r);
   \end{equation*}
  \item $\mass{T_j}\leq (1+\lambda)\nu\gamma\mass{T}$ for all $j\in\N$, where $\nu:=\delta$ if $k=1$ or $\nu:=\max\{\delta^k,\delta^\alpha\}$ otherwise;
  \item $\diam(\spt T_j)\leq 4G\left(\gamma^{-1}\frac{2}{1-\lambda}5^{k+\alpha}\mass{T_j}\right)$;
  \item $\mass{R}+\frac{1-\lambda}{1+\lambda}\sum_{i=1}^\infty\mass{T_i}\leq \mass{T}$.
 \end{enumerate} 
\et
For the exact value of $\gamma$ see the beginning of the proof of \propref{Proposition:decomposing-with-omega}. If $k=1$, all statements of the theorem hold for $\lambda=0$ as well.
In the proof of \thmref{theorem:suitable-decomposition} we will need the following result, 
which generalizes Lemma 3.4 of \cite{Wenger-GAFA} and partially Lemma 3.1 of \cite{Wenger-filling}.  

\bp\label{proposition:isop-ineq-growth-estimate}
 Let $X$ be a complete metric space, $k\geq 2$, $\alpha>1$, and suppose that $X$ admits an isoperimetric inequality of rank 
 $\alpha$ for $\intcurr_{k-1}(X)$ with a constant $D_{k-1}\in [1,\infty)$.
 Then for every $T\in\intcurr_{k-1}(X)$ with $\bdry T=0$ and every $\varepsilon>0$ there exists an $S\in\intcurr_k(X)$ with $\bdry S=T$, satisfying
 \begin{equation}\label{equation:upper-bound-mass-volgrowth-estimate-prop}
  \mass{S}\leq \min\left\{(1+\varepsilon)\fillvol(T), D_{k-1}I_{k,\alpha}(\mass{T})\right\}
 \end{equation}
 and with the following property: For every $x\in \spt S$ and every $0\leq r\leq \dist(x, \spt T)$ we have
  \begin{equation*}
  \|S\|(B(x,r))\geq F_{k,\alpha,\mu}(r)
 \end{equation*}
 where 
 \begin{equation*}
  \mu:= \min\left\{\frac{1}{(3D_{k-1})^{k-1}\alpha_1^k}, \frac{1}{(3D_{k-1})^{\alpha-1}\alpha_1^\alpha}\right\}
 \end{equation*}
 with $\alpha_1:=\max\{k,\alpha\}$.
\ep

The proof relies on the arguments contained in \cite[Theorem 10.6]{Ambr-Kirch-curr}.

\begin{proof}
Let ${\mathcal M}$ denote the complete metric space consisting of all $S\in\intcurr_k(X)$ with $\bdry S=T$ and endowed with the 
metric given by $d_{\mathcal M}(S,S'):= \mass{S-S'}$. 
Choose an $\tilde{S}\in{\mathcal M}$ satisfying
\eqref{equation:upper-bound-mass-volgrowth-estimate-prop}. By a well-known variational principle (see e.g.~\cite{Ekeland}) there exists an $S\in{\mathcal M}$ with 
$\mass{S}\leq \mass{\tilde{S}}$ 
and such that the function 
\begin{equation*}
 S'\mapsto \mass{S'}+\frac{1}{2}\mass{S-S'}
\end{equation*}
has a minimum at $S'=S$. Let $x\in\spt S\ohne\spt T$ and set $R:= \dist(x,\spt\bdry T)$. 
We claim that if $r\in(0,R)$ then
  \begin{equation}\label{equation:explicit-isop-growth-estimate}
  \|S\|(B(x,r))\geq \left\{\begin{array}{l@{\quad}l}
     \frac{r^k}{(3D_{k-1})^{k-1}\alpha_1^k} & r\leq 3D_{k-1}\alpha_1\\
     \frac{r^\alpha}{(3D_{k-1})^{\alpha-1}\alpha_1^\alpha} & r> 3D_{k-1}\alpha_1.
   \end{array}\right.
 \end{equation}
First note that the slicing theorem \cite[Theorems 5.6 and 5.7]{Ambr-Kirch-curr} implies that for almost every $r\in(0,R)$
the slice $\bdry(S\rstr B(x,r))$ exists, has zero boundary, and belongs to $\intcurr_{k-1}(X)$. For an
$S_r\in\intcurr_k(X)$ with $\bdry S_r=\bdry(S\rstr B(x,r))$ the integral current
$S\rstr (X\ohne B(x,r))+S_r$ has boundary $T$ and thus, comparison with $S$ yields
\begin{equation*}
\mass{S\rstr (X\ohne B(x,r))+S_r} + \frac{1}{2}\mass{S\rstr B(x,r)-S_r}\geq \mass{S}.
\end{equation*}
If, moreover, $S_r$ is chosen such that $\mass{S_r}\leq D_{k-1}I_{k,\alpha}(\mass{\bdry(S\rstr B(x,r))})$ then it follows that
\begin{equation*}
\mass{S\rstr B(x,r)}\leq3\mass{S_r}\leq 3D_{k-1}I_{k,\alpha}(\mass{\bdry(S\rstr B(x,r))})                                                   
\end{equation*}
and consequently, 
\begin{equation}\label{equation:diff-inequality-beta}
 \beta(r)\leq 3D_{k-1}I_{k,\alpha}(\beta'(r))
\end{equation}
for almost every $r\in(0,R)$, where $\beta(r):= \|S\|(B(x,r))$.

Set $\overline{r}:=\sup\{r\in[0,R]: \beta(r)\leq 3D_{k-1}\}$ and observe that for almost every $r\in(0,\overline{r})$
\begin{equation*}
 \frac{d}{dr}\left[\beta(r)^{\frac{1}{k}}\right]=\frac{\beta'(r)}{k\,\beta(r)^{\frac{k-1}{k}}}\geq \frac{1}{(3D_{k-1})^{\frac{k-1}{k}}k}.
\end{equation*}
This yields
\begin{equation*}
 \beta(r)\geq \frac{r^k}{(3D_{k-1})^{k-1}k^k}
\end{equation*}
for all $r\in[0,\overline{r}]$ and consequently
\begin{equation*}
 \beta(r)\geq \frac{r^k}{(3D_{k-1})^{k-1}\alpha_1^k}
\end{equation*}
for all $r\in[0,\overline{R}]$, where $\overline{R}:=\min\{R,3D_{k-1}\alpha_1\}$.
Indeed, it is clear that $\overline{r}\leq \overline{R}$ and in case $\overline{r}<\overline{R}$ we furthermore have 
\begin{equation*}
 \beta(r)\geq \beta(\overline{r})\geq 3D_{k-1}\geq \frac{r^k}{(3D_{k-1})^{k-1}\alpha_1^k}
\end{equation*}
for all $r\in[\overline{r}, \overline{R}]$.
This proves \eqref{equation:explicit-isop-growth-estimate} for $r\in[0,\overline{R}]$. Now, if $R>3D_{k-1}\alpha_1$ then for almost every
$r\in [3D_{k-1}\alpha_1, R]$
\begin{equation*}
 3D_{k-1}\leq \beta(\overline{r})\leq \beta(r)\leq 3D_{k-1}I_{k,\alpha}(\beta'(r))
\end{equation*}
and hence $\beta'(r)\geq 1$. It follows that
\begin{equation*}
 \frac{d}{dr}\left[\beta(r)^{\frac{1}{\alpha}}\right]=\frac{\beta'(r)}{\alpha\,\beta(r)^{\frac{\alpha-1}{\alpha}}}
   \geq \frac{1}{(3D_{k-1})^{\frac{\alpha-1}{\alpha}}\alpha_1}
\end{equation*}
and thus
\begin{equation*}
 \beta(r)^{\frac{1}{\alpha}}\geq \beta(3D_{k-1}\alpha_1)^{\frac{1}{\alpha}} + \frac{r-3D_{k-1}\alpha_1}{(3D_{k-1})^{\frac{\alpha-1}{\alpha}}\alpha_1}
  \geq \frac{r}{(3D_{k-1})^{\frac{\alpha-1}{\alpha}}\alpha_1}.
\end{equation*}
This concludes the proof of \eqref{equation:explicit-isop-growth-estimate}. In order to finish the proof of the proposition it is enough to 
show the statement for $r\in[1,3D_{k-1}\alpha_1]$, since the other cases are direct consequences of \eqref{equation:explicit-isop-growth-estimate}.
We simply calculate
\begin{equation*}
 \|S\|(B(x,r))\geq \frac{r^k}{(3D_{k-1})^{k-1}\alpha_1^k} \geq \frac{r^\alpha}{(3D_{k-1})^{k-1}\alpha_1^k(3D_{k-1}\alpha_1)^{\alpha-k}}
\end{equation*}
to obtain the desired inequality.
\end{proof}

A direct consequence of the proposition is the following estimate on the filling radius.

\bc\label{corollary:fillrad-fillvol}
 Let $X$ be a complete metric space, $k\geq 2$, $\alpha>1$, and suppose that $X$ admits an isoperimetric inequality of rank $\alpha$ for $\intcurr_{k-1}(X)$. 
 Then for every $T\in\intcurr_{k-1}(X)$ with $\bdry T=0$ we have
 \begin{equation*}
  \fillrad_X(T)\leq G_{k,\alpha}(\mu^{-1}\fillvol_X(T))\leq\left\{\begin{array}{l@{\qquad}l}
   \mu'\mass{T}^{\frac{1}{k-1}} &\mass{T}\leq 1\\
   \mu'\mass{T}^{\frac{1}{\alpha-1}} &\mass{T}>1,
  \end{array}\right.
 \end{equation*}
 where 
 \begin{equation*}
  \mu':=\max\left\{\left(\frac{D_k}{\mu}\right)^{\frac{1}{k}}, \left(\frac{D_k}{\mu}\right)^{\frac{1}{\alpha}}\right\}.
 \end{equation*}
\ec

\subsection{An analytic lemma} 
For $(k,\alpha)\in\Lambda$ and $\gamma\in(0,\infty)$ we first define
an auxiliary function by
\begin{equation*}
 H_{1,0,\gamma}(r)= \gamma 
\end{equation*}
and
\begin{equation*}
 H_{k,\alpha,\gamma}(r):=\left\{\begin{array}{l@{\qquad}l}
  \gamma^{\frac{1}{k}}\cdot r^{\frac{k-1}{k}} & 0\leq r\leq \gamma\\
  \gamma^{\frac{1}{\alpha}}\cdot r^{\frac{\alpha-1}{\alpha}} & \gamma < r <\infty
 \end{array}\right.
\end{equation*}
if $k\geq 2$.
For the convenience of the reader we summarize some simple properties of the auxiliary functions thus far defined. Their properties 
will be used in the sequel, sometimes without explicit mentioning.
\bl\label{lemma:trivial-properties-functions}
 Let $(k,\alpha)\in\Lambda$ and $\gamma\in(0,\infty)$ and set $F:=F_{k,\alpha,\gamma}$, $G:=G_{k,\alpha}$, $H:=H_{k,\alpha,\gamma}$, and $I:= I_{k,\alpha}$. 
 Then the following properties hold:
 \begin{enumerate}
  \item For all $r\geq 0$ we have $F(5r)\leq 5^{k+\alpha}F(r)$;
  \item If $k\geq 2$ and $\nu\geq 0$ then
   \begin{equation*}
    \min\left\{\nu^{\frac{1}{k}}, \nu^{\frac{1}{\alpha}}\right\} G(r)\leq G(\nu r)\leq \max\left\{\nu^{\frac{1}{k}}, \nu^{\frac{1}{\alpha}}\right\} G(r)
   \end{equation*}
   for all $r\geq 0$;
  \item If $k\geq 2$ and $\nu\geq 0$ then 
   \begin{equation*}
    \gamma\min\left\{\nu^k,\nu^\alpha\right\}\;r\leq F(\nu G(r))\leq \gamma\max\left\{\nu^k, \nu^\alpha\right\}\; r
   \end{equation*}
   for all $r\geq 0$;
  \item If $k\geq 2$ and $\nu\geq 0$ then 
   \begin{equation*}
    \min\left\{(\gamma\nu)^{\frac{1}{k}},(\gamma\nu)^{\frac{1}{\alpha}}\right\}\;r\leq G(\nu F(r))
     \leq \max\left\{(\gamma\nu)^{\frac{1}{k}}, (\gamma\nu)^{\frac{1}{\alpha}}\right\}\; r
   \end{equation*}
   for all $r\geq 0$;
  \item If $k=1$ then $F'(r)= H(F(r))$ for all $r\geq 0$
  \item If $k\geq 2$ then $F'(r)=kH(F(r))$ when $r\in (0,1)$ and $F'(r)=\alpha H(F(r))$ when $r>1$;
  \item If $k\geq 2$ then $I(s) + I(t)\leq I(s+t)$ for all $s,t\geq 0$;
  \item If $k\geq 2$ and $1\leq \mu\leq 1/\gamma$ then
   \begin{equation*}
    I(\mu H(r)) \leq \mu \max\left\{(\mu\gamma)^{\frac{1}{k-1}}, (\mu\gamma)^{\frac{1}{\alpha-1}}\right\}\cdot r
   \end{equation*}
   for all $r\geq 0$.
 \end{enumerate}
\el

The proof is by straight-forward verification and is therefore omitted. For the proof of \thmref{theorem:suitable-decomposition} we need 
the following analytic lemma.
\bl\label{Lemma:slow-growth}
 Let $(k,\alpha)\in\Lambda$, $\gamma\in(0,1)$, and abbreviate $F:= F_{k,\alpha,\gamma}$ and $H:=H_{k,\alpha,\gamma}$. Let furthermore $r_0>0$ and suppose 
 $f:[0,r_0]\to [0,\infty)$ is non-decreasing and continuous from the right with 
 $f(r_0)<5^{-(k+\alpha)}F(r_0)$ 
and such that
 \begin{equation*}
  r_*:= \max\left\{r\in[0,r_0]: f(r)\geq F(r)\right\} >0.
 \end{equation*}
 Then $r_*< r_0/5$ and there is a measurable subset $K\subset (r_*,r_0/5)$ of strictly positive Lebesgue measure such that
 \begin{equation*}
  f(5r)<5^{k+\alpha} f(r)\quad\text{ and }\quad f'(r)< (k+\alpha)H(f(r))
 \end{equation*}
 for every $r\in K$.
\el

This lemma will be applied with $f(r)$ the mass in a ball of radius $r$ of an integral current.
\begin{proof}
 First of all, if $r_*\geq r_0/5$ then it follows that
 \begin{equation*}
  F(r_0)\leq F(5r_*)\leq 5^{k+\alpha}F(r_*)=5^{k+\alpha}f(r_*)\leq 5^{k+\alpha}f(r_0),
 \end{equation*} 
 which contradicts the hypothesis. This proves that indeed $r_*<r_0/5$.
 Now suppose that for almost every $r\in (r_*, r_0/5)$ we have 
 \begin{equation*}
  \text{either}\quad f(5r)\geq 5^{k+\alpha}f(r)\quad\text{or}\quad f'(r)\geq (k+\alpha)H(f(r)).
 \end{equation*}
Define
 \begin{equation*}
  r'_*:= \inf\left\{r\in[r_*, r_0/5]: f(5r)\geq 5^{k+\alpha}f(r)\right\},
 \end{equation*}
where we agree on $\inf\emptyset = \infty$.
 It then follows that $r'_*>r_*$ since otherwise
 \begin{equation*}
 F(5r_*)\leq 5^{k+\alpha}F(r_*)=5^{k+\alpha}f(r_*)\leq f(5r_*),
 \end{equation*}
 in contradiction with the definition of $r_*$.
 If $k=1$ then set $r''_*:=\min\{r'_*,r_0/5\}$ and note that $f'(r)\geq \gamma$ for almost every $r\in(r_*,r''_*)$ and thus 
 $f(r''_*)\geq f(r_*)+\gamma(r''_*-r_*)=\gamma r''_*$, which is impossible.
 If, on the other hand, $k\geq 2$ then we distinguish the following two cases. 
 Suppose first that $r_*<1$ and set $r''_*:= \min\{1, r_0/5, r'_*\}$; observe that $r''_*>r_*$ and $f(r''_*)< \gamma$. 
 Consequently, we have 
 \begin{equation*}
  \frac{d}{dr}\left[ f(r)^{\frac{1}{k}}\right] = \frac{f'(r)}{kf(r)^{\frac{k-1}{k}}} \geq \frac{(k+\alpha)H(f(r))}{kf(r)^{\frac{k-1}{k}}} > \gamma^{\frac{1}{k}}
 \end{equation*}
 for almost every $r\in(r_*,r''_*)$ and hence
 \begin{equation*}
  f(r''_*)^{\frac{1}{k}}> f(r_*)^{\frac{1}{k}} + \gamma^{\frac{1}{k}}(r''_*-r_*)=\gamma^{\frac{1}{k}}r''_*,
 \end{equation*}
 which is not possible. 
 Suppose next that $r_*\geq 1$ and set $r''_*:= \min\{r_0/5, r'_*\}$; 
 observe that $r''_*>r_*$ and $f(r''_*)>\gamma$, from which we conclude analogously as above
 that 
 \begin{equation*}
  \frac{d}{dr}\left[ f(r)^{\frac{1}{\alpha}}\right] = \frac{f'(r)}{\alpha f(r)^{\frac{\alpha-1}{\alpha}}} 
          \geq \frac{(k+\alpha)H(f(r))}{\alpha f(r)^{\frac{\alpha-1}{\alpha}}} > \gamma^{\frac{1}{\alpha}}
 \end{equation*}
for almost every $r\in(r_*,r_*'')$ and thus
 \begin{equation*}
  f(r''_*)^{\frac{1}{\alpha}}> f(r_*)^{\frac{1}{\alpha}} + \gamma^{\frac{1}{\alpha}}(r''_*-r_*)>\gamma^{\frac{1}{\alpha}}r''_*,
 \end{equation*}
 again a contradiction with the definition of $r_*$. This concludes the proof of the lemma.
\end{proof}

\subsection{Controlling the thin parts of a current}
Let $X$ be a complete metric space and fix $(k,\alpha)\in\Lambda$.
The following set which we associate with an element $T\in\intcurr_k(X)$ and constants $\gamma\in(0,1)$ and $L\in(0,\infty]$ will sometimes be referred to as 
the thin part of $T$,
\begin{equation*}
\begin{split}
 \Omega(T,\gamma, L):= \Big\{ x\in\spt T: \;\Theta_{*k}(\|T\|, x)>\frac{\gamma}{\omega_k}&\text{ and }
  \|T\|(B(x,r))<\frac{1}{2}5^{-(k+\alpha)}F_{k,\alpha,\gamma}(r)\\
 &\text{ for an $r\in\left[0,\min\{L, \dist(x,\spt\bdry T)\}\right]$}\Big\}.
\end{split}
\end{equation*}
Here, $\omega_k$ denotes the volume of the unit ball in $\R^k$. Note that in the above we explicitly allow the value $L=\infty$. Furthermore, we agree on the convention $\dist(x,\emptyset)=\infty$. It should be remarked
that $\Omega(T,\gamma,L)$ also depends on $\alpha$ even though we omit $\alpha$ in our notation. 
The inequality involving the lower
density is satisfied for $\|T\|$-almost every $x\in\spt T$ if $\gamma<\omega_kk^{-k/2}$ by \cite{Ambr-Kirch-curr}. It is not difficult to see that
$\Omega(T,\gamma,L)$ is then $\|T\|$-measurable and that, in case $\bdry T=0$, we have $\Omega(T,\gamma,\infty)=\spt T$ up to a set of $\|T\|$-measure zero.

\bl\label{Lemma:suitable-points-on-support}
 Let $X$ be a complete metric space, $(k,\alpha)\in\Lambda$, and $\gamma\in(0,a)$, where $a:=\min\{1, \omega_kk^{-k/2}\}$. 
 Abbreviate $F:=F_{k,\alpha,\gamma}$, $G:= G_{k,\alpha}$ and $H:=H_{k,\alpha,\gamma}$.
 Let furthermore $T\in\intcurr_k(X)$ 
 and $L\in(0,\infty]$. Then there exist finitely many points $x_1,\dots, x_N\in\Omega(T,\gamma,L)$ and
 $s_1,\dots, s_N\in (0, \infty)$ with the following properties:
 \begin{enumerate}
  \item With $A:= G(\gamma^{-1}\|T\|(B(x_i,s_i)))$ we have
    \begin{equation*}
       A< s_i<\min\left\{\frac{L}{5}, \frac{1}{5}\dist(x_i, \spt\bdry T), 2\cdot 5^{k+\alpha}A\right\}
    \end{equation*}
  \item $B(x_i, 2s_i)\cap B(x_j, 2s_j)=\emptyset$ for all $i\not=j$
  \item $T\rstr B(x_i, s_i)\in\intcurr_k(X)$
  \item $\frac{1}{2}5^{-(k+\alpha)}F(s_i)\leq\|T\|(B(x_i, s_i))\leq F(s_i)$
  \item $\mass{\bdry(T\rstr B(x_i, s_i))}\leq (k+\alpha)H(\|T\|(B(x_i, s_i)))$
  \item $\sum_{i=1}^N\|T\|(B(x_i, s_i))\geq 5^{-(k+\alpha)}\|T\|(\Omega(T,\gamma,L))$.
 \end{enumerate}
\el

We note that in the above we allow $N=0$ if $\|T\|(\Omega(T,\gamma,L))=0$.

\begin{proof}
 For each $x\in\Omega(T,\gamma, L)$ set
 \begin{equation*}
  f_x(r):= \|T\|(B(x,r))\qquad\text{for $r\in[0, \infty)$}
 \end{equation*}
and note that $f_x$ is non-decreasing and continuous from the right.
Define furthermore
 \begin{equation*}
  r_0(x):= \inf\left\{ r\in[0,\min\{L, \dist(x,\spt\bdry T)\}]: \|T\|(B(x,r))< \frac{1}{2}5^{-(k+\alpha)}F(r)\right\}.
 \end{equation*}
Since
 \begin{equation*}
  \liminf_{r\searrow 0}\frac{f_x(r)}{r^k} = \omega_k\Theta_{*k}(\|T\|,x)>\gamma,
 \end{equation*}
 it follows that $r_0(x)>0$ and
 \begin{equation*}
  r_*(x):= \max\left\{r\in[0, r_0(x)]: f_x(r)\geq F(r)\right\}>0.
 \end{equation*}
 Note that we also have
 \begin{equation}\label{equation:bound-r0}
  r_0(x)=G(\gamma^{-1}F(r_0(x)))= G(2\gamma^{-1}5^{k+\alpha}\mass{T})
 \end{equation}
 and $f_x(r_0(x))< 5^{-(k+\alpha)}F(r_0(x))$.
 \lemref{Lemma:slow-growth} and the slicing theorem for rectifiable currents imply that there exists
 for each $x\in\Omega(T,\gamma,L)$ an $r(x)\in(r_*(x), r_0(x)/5)$ such that
 \begin{enumerate}
  \item[(a)] $T\rstr B(x, r(x))\in\intcurr_k(X)$
  \item[(b)] $\|T\|(B(x, r(x)))< F(r(x))$
  \item[(c)] $\|T\|(B(x, 5r(x)))<5^{k+\alpha}\|T\|(B(x,r(x)))$
  \item[(d)] $\mass{\bdry(T\rstr B(x, r(x)))}\leq f'_x(r(x))<(k+\alpha)H(\|T\|(B(x,r(x))))$.
 \end{enumerate}
 The points $x_1,\dots, x_N$ and the radii $s_1,\dots,s_N$ are now constructed as follows: Set $\Omega_1:= \Omega(T,\gamma,L)$ and 
 $s_1^*:= \sup\{ r(x): x\in\Omega_1\}$. From \eqref{equation:bound-r0} it follows that $s^*_1<\infty$.
Choose $x_1\in\Omega_1$ in such a way that $r(x_1)>\frac{2}{3}s_1^*$. If $x_1,\dots, x_j$ are chosen define
\begin{equation*}
  \Omega_{j+1}:= \Omega(T,\gamma,L)\ohne\bigcup_{i=1}^jB(x_i, 5r(x_i))
\end{equation*}
and
\begin{equation*}
 s_{j+1}^*:= \sup\{r(x): x\in\Omega_{j+1}\}.
\end{equation*}
If $\|T\|(\Omega_{j+1})>0$ we can choose $x_{j+1}\in\Omega_{j+1}$ such that $r(x_{j+1})>\frac{2}{3}s_{j+1}^*$. This procedure yields (possibly finite)
sequences $x_j\in\Omega_j$, $s_1^*\geq s_2^*\geq \dots\geq 0$, and $s_i:= r(x_i)$. 
We show that for a suitably large $N$ the so defined points and numbers have the desired properties stated in the lemma. 
We first note that, by (b) and the definition of $r_0(x)$,
\begin{equation*}
 \frac{1}{2}5^{-(k+\alpha)}F(s_i)\leq \|T\|(B(x_i,s_i))<F(s_i),
\end{equation*}
which proves (iv). Property (i) follows from this and the fact that $s_i=G(\gamma^{-1}F(s_i))$.
Furthermore, we have 
\begin{equation*}
 d(x_i, x_{i+\ell})> 5s_i = 2s_i + 3s_i > 2s_i+ 2s_i^*\geq 2s_i+ 2s_{i+\ell}
\end{equation*}
and thus we obtain (ii). Properties (iii) and (v) are direct consequences of (a) and (d), respectively. 
We are therefore left to show
that (vi) holds for some $N\in\N$. On the one hand, if $\|T\|(\Omega_{n+1})=0$ for some $n\in\N$ then (c) yields
\begin{equation*}
 \sum_{i=1}^n\|T\|(B(x_i, s_i))> 5^{-(k+\alpha)}\sum_{i=1}^n\|T\|(B(x_i, 5s_i))\geq 5^{-(k+\alpha)}\|T\|(\Omega(T,\gamma,L)),
\end{equation*}
which establishes (vi) and thus the lemma with $N=n$. On the other hand, if $\|T\|(\Omega_n)>0$ for all $n\in\N$ then it follows easily that $s_n^*\searrow 0$. 
Indeed, this is a consequence of the fact that
\begin{equation*}
 \frac{1}{2}5^{-(k+\alpha)}\sum_{i=1}^\infty F\left(\frac{2}{3}s_i^*\right)<\sum_{i=1}^\infty\frac{1}{2}5^{-(k+\alpha)}F(s_i)
 \leq\sum_{i=1}^\infty \|T\|(B(x_i, s_i))\leq \mass{T}<\infty.
\end{equation*}
Furthermore we claim that 
\begin{equation*}
 \|T\|\left(\Omega(T,\gamma,L)\ohne\bigcup_{i=1}^\infty B(x_i, 5s_i)\right)=0.
\end{equation*}
If this were not true we would have $x\in\Omega(T,\gamma,L)\ohne\cup_{i=1}^\infty B(x_i, 5s_i)$ and since $r(x)>0$ we would obtain a contradiction with 
$s_i^*\searrow 0$. The rest now follows as in the case above.
\end{proof}

\subsection{Proof of \thmref{theorem:suitable-decomposition}}
The following proposition is an intermediate step on the way to the proof of the main theorem of this section. The proposition shows how to
construct a suitable decomposition of a current $T\in\intcurr_k(X)$ in a way that helps to reduce the set $\Omega(T,\gamma, L)$.
\bp\label{Proposition:decomposing-with-omega}
 Let $X$ be a complete metric, $(k,\alpha)\in\Lambda$ and suppose in case $k\geq 2$ that 
 $X$ admits an isoperimetric inequality of rank $\alpha$ for $\intcurr_{k-1}(X)$.  
 For every $\lambda\in(0,1)$ there exists $\gamma\in(0,1)$ with the following property. Set $F:=F_{k,\alpha,\gamma}$ and $G:=G_{k,\alpha}$.
 Then for every $L\in(0,\infty]$ and $T\in\intcurr_k(X)$ there is a decomposition
 \begin{equation*}
  T = R+T_1+\dots+T_N
 \end{equation*}
 with $R, T_i\in\intcurr_k(X)$ such that
 \begin{enumerate}
  \item $\bdry R=\bdry T$ and $\bdry T_i=0$
  \item $\mass{T_i}\leq (1+\lambda)F(L/5)$
  \item $\diam(\spt T_i)\leq 4G\left(\gamma^{-1}\frac{2}{1-\lambda}5^{k+\alpha}\mass{T_i}\right)$
  \item $\mass{R}+ \frac{1-\lambda}{1+\lambda}\sum_{i=1}^N\mass{T_i}\leq \mass{T}$
  \item $\sum_{i=1}^N\mass{T_i}\geq (1-\lambda)5^{-(k+\alpha)}\|T\|(\Omega(T,\gamma,L))$.
 \end{enumerate}
\ep

We will prove that for $k=1$ all the properties above also hold with $\lambda=0$.
\begin{proof}
 If $k=1$ then set $\gamma:=1/2$. If $k\geq 2$ then define
 \begin{equation*}
  \gamma:= \frac{1}{k+\alpha}\min\left\{A^{k-1}, A^{\alpha-1}, \omega_kk^{-\frac{k}{2}}\right\}\quad\text{ with }\quad A:=\frac{\lambda}{3D_{k-1}(k+\alpha)}
 \end{equation*}
 and where $D_{k-1}$ denotes the constant in the isoperimetric inequality for $\intcurr_{k-1}(X)$. We may of course assume that 
 $D_{k-1}\geq 1$. 
 We may furthermore assume that $\|T\|(\Omega(T,\gamma,L))>0$ since otherwise we can set $R:=T$ and there is then 
 nothing to prove. Let $x_1,\dots,x_N\in\Omega(T,\gamma,L)$ and $s_1,\dots,s_N\in(0,\infty)$ be as in 
 \lemref{Lemma:suitable-points-on-support}. Fix $i\in\{1,\dots, N\}$. 
 If $k=1$ then set $T_i:= T\rstr B(x_i,s_i)$ and note that $\mass{\bdry T_i}\leq \gamma<1$ and thus $\bdry T_i=0$.
 If, on the other hand, $k\geq 2$ then
 choose $S_i\in\intcurr_k(X)$ such that $\bdry S_i= \bdry(T\rstr B(x_i,s_i))$ and with the properties of 
 \propref{proposition:isop-ineq-growth-estimate}. We have 
 \begin{equation}\label{equation:mass-cup}
  \mass{S_i}\leq D_{k-1}I_{k,\alpha}(\mass{\bdry(T\rstr B(x_i,s_i))})\leq \lambda\|T\|(B(x_i,s_i)),
 \end{equation}
where the second inequality follows from (viii) of \lemref{lemma:trivial-properties-functions} and the definition of $\gamma$. 
 Next we have that $\spt S_i\subset B(x_i, 2s_i)$. This is indeed a consequence of \propref{proposition:isop-ineq-growth-estimate}, (iv) of \lemref{lemma:trivial-properties-functions}, the fact
 that
 \begin{equation*}
  \mass{S_i}\leq \lambda \|T\|(B(x_i,s_i))\leq\lambda F(s_i)=\left\{\begin{array}{l@{\quad}l}
     \lambda\gamma\cdot s_i^k &  s_i\leq 1\\
     \lambda\gamma\cdot s_i^\alpha & s_i>1,
   \end{array}\right.
 \end{equation*}
and the choice of $\gamma$.
Thus $T_i:= T\rstr B(x_i, s_i) - S_i$ satisfies $T_i\in\intcurr_k(X)$, $\bdry T_i= 0$ and $\spt T_i\subset B(x_i, 2s_i)$. 
From \eqref{equation:mass-cup} we see that
 \begin{equation}\label{equation:mass-cycles}
  (1-\lambda)\|T\|(B(x_i, s_i))\leq \mass{T_i}\leq (1+\lambda)\|T\|(B(x_i, s_i))
 \end{equation}
 and thus
 \begin{equation*}
  \mass{T_i}\leq (1+\lambda)\|T\|(B(x_i, s_i))\leq (1+\lambda)F(s_i)\leq (1+\lambda)F(L/5),
 \end{equation*}
 which proves (ii) of the present proposition. 
 Note that the above conclusion holds with $\lambda=0$ in the case $k=1$.
 We proceed as above for every $i\in\{1,\dots, N\}$ and note that in each step of the construction only the ball $B(x_i, 2s_i)$, which
 is disjoint from the other balls, is affected. We thus obtain cycles $T_1,\dots, T_N$ and we claim that these together with $R:= T-T_1 -\dots - T_N$ have the 
 properties stated in the proposition.
 Indeed, (i) is obvious and (ii) has already been proved. As for (iii) it is enough to note that $\diam(\spt T_i)\leq 4s_i$ and that, by (iv) of \lemref{lemma:trivial-properties-functions} and by \eqref{equation:mass-cycles},
 \begin{equation*}
  s_i = G(\gamma^{-1}F(s_i))\leq G\left(\gamma^{-1}\frac{2}{1-\lambda} 5^{k+\alpha}\mass{T_i}\right).
 \end{equation*}
 Again, if $k=1$ then the above holds with $\lambda=0$.
 Now, our construction yields
 \begin{equation*}
  \begin{split}
   \mass{R}&\leq \|T\|\left(X\ohne\bigcup_{i=1}^NB(x_i, s_i)\right) + \lambda\sum_{i=1}^N\|T\|(B(x_i, s_i))\\
           &= \mass{T} -(1-\lambda)\sum_{i=1}^N\|T\|(B(x_i, s_i))\\
           &\leq \mass{T} - \frac{1-\lambda}{1+\lambda}\sum_{i=1}^N\mass{T_i}
  \end{split}
 \end{equation*}
 from which (iv) follows. Note that in the case $k=1$ we have $R=T\rstr(X\backslash\cup_{i=1}^N B(x_i,s_i))$ and $T_i=T\rstr B(x_i, s_i)$ and hence
 \begin{equation*}
  \mass{R} + \sum_{i=1}^N\mass{T_i} = \|T\|\left(X\ohne\bigcup_{i=1}^NB(x_i, s_i)\right) + \sum_{i=1}^N\|T\|(B(x_i, s_i))=\mass{T}.
 \end{equation*}
 Finally, we use (vi) of \lemref{Lemma:suitable-points-on-support} together with \eqref{equation:mass-cycles} to calculate
 \begin{equation*}
  \sum_{i=1}^N\mass{T_i}\geq (1-\lambda)\sum_{i=1}^N\|T\|(B(x_i, s_i))\geq (1-\lambda)5^{-(k-\alpha)}\|T\|(\Omega(T,\gamma,L)).
 \end{equation*}
 This establishes (v) and concludes the proof of the proposition.
\end{proof}

We are now ready for the proof of the decomposition theorem stated at the beginning of the section.
\begin{proof}[{Proof of \thmref{theorem:suitable-decomposition}}]
 Let $\gamma$ be as in \propref{Proposition:decomposing-with-omega}.
 Set $R_0:= T$ and $N_0:=0$. Successive application of \propref{Proposition:decomposing-with-omega} yields possibly finite sequences
 $(R_i), (T_j)\subset\intcurr_k(X)$ and a strictly increasing sequence of integers $N_1<N_2<\dots$ such that for every $i\in\N\cup\{0\}$
 \begin{equation*}
  R_i= R_{i+1}+ T_{N_i+1}+\dots+T_{N_{i+1}}
 \end{equation*}
 and such that the following properties hold:
 \begin{enumerate}
  \item[(a)] $\bdry R_i=\bdry T$ and $\bdry T_j=0$ for all $i, j$
  \item[(b)] $\mass{T_j}\leq (1+\lambda)\nu\gamma\mass{R_i}$ for all $j\in\{N_i+1,\dots,N_{i+1}\}$
  \item[(c)] $\diam(\spt T_j)\leq 4G\left(\gamma^{-1}\frac{2}{1-\lambda}5^{k+\alpha}\mass{T_j}\right)$
  \item[(d)] $\mass{R_{i+1}}+\frac{1-\lambda}{1+\lambda}\sum_{j=N_i+1}^{N_{i+1}}\mass{T_j}\leq \mass{R_i}$
  \item[(e)] $\sum_{j=N_i+1}^{N_{i+1}}\mass{T_j}\geq (1-\lambda)5^{-(k+\alpha)}\|R_i\|(\Omega(R_i,\gamma,L_i))$.
 \end{enumerate}
Here, $L_i$ is defined by $L_i:= 5\delta G(\mass{R_i})$ and $\nu:=\delta$ if $k=1$ and $\nu=\max\{\delta^k,\delta^\alpha\}$ otherwise.
Property (b) follows from (iii) of \lemref{lemma:trivial-properties-functions} and the definition of $L_i$. We note that in the case $k=1$, all the properties above hold with $\lambda=0$.
We thus obtain for each $i\in\N\cup\{0\}$ a decomposition 
\begin{equation*}
 T=R_i + \sum_{j=1}^{N_i}T_j
\end{equation*}
which, by property (d), satisfies
\begin{equation}\label{equation:mass-decomp-estimate}
 \mass{R_i}+ \frac{1-\lambda}{1+\lambda}\sum_{j=1}^{N_i}\mass{T_j}\leq \mass{T},
\end{equation}
or, if $k=1$,
\begin{equation*}
 \mass{R_i} + \sum_{j=1}^{N_i}\mass{T_j} = \mass{T}.
\end{equation*}
In particular, we have
\begin{equation*}
 \mass{R_{i+m} - R_i}= \mass{T_{N_i+1}+\dots+T_{N_{i+m}}}\leq \sum_{j=N_i+1}^\infty\mass{T_j}\to 0
\end{equation*}
as $i\to\infty$, thus the sequence $(R_i)$ is Cauchy with respect to the mass norm. Since the additive group of integer rectifiable 
$k$-currents together with the mass norm is complete, there exists 
$R\in\intcurr_k(X)$ such that $\mass{R-R_i}\to0$ and, in particular,
\begin{equation*}
 T = R + \sum_{j=1}^\infty T_j.
\end{equation*}
Clearly, we have $\bdry R=\bdry T$ and thus property (i) holds. 
Properties (iii), (iv) and (v) are direct consequences of (b), (c) and \eqref{equation:mass-decomp-estimate}, respectively. We are therefore left to 
establish (ii).
For this let $x\in\spt R\ohne\spt\bdry T$ and $$0<r<\min\{5\delta G(\mass{R}), \dist(x,\spt\bdry T)\}.$$ 
Observe that 
\begin{equation*}
 \|R_i\|(B(x,t))\to \|R\|(B(x,t))
\end{equation*}
and $\|R\|(B(x,t))>0$ for all $t\in(0,r]$.
Fix $0<s<r$ and $\varepsilon>0$. By (e) and \eqref{equation:mass-decomp-estimate} we have 
 \begin{equation*}
  \|R_i\|(\Omega(R_i,\gamma, L_i))\to 0.
 \end{equation*}
There thus exist $i_0\in\N$ and $x'\in\spt R_{i_0}$ with $d(x,x')\leq s$, 
\begin{equation*}
\|R\|(B(x,r))\geq (1-\varepsilon)\|R_{i_0}\|(B(x,r))
\end{equation*}
and such that
\begin{equation*}
 \|R_{i_0}\|(B(x',r-s))\geq \frac{1}{2}5^{-(k+\alpha)}F(r-s).
\end{equation*}
It finally follows that
\begin{equation*}
 \begin{split}
  \|R\|(B(x,r))&\geq (1-\varepsilon)\|R_{i_0}\|(B(x,r))\\
               &\geq (1-\varepsilon)\|R_{i_0}\|(B(x', r-s))\\
               &\geq \frac{1}{2}(1-\varepsilon)5^{-(k+\alpha)}F(r-s).
 \end{split}
\end{equation*}
Since $s$ and $\varepsilon$ were arbitrary this establishes (ii) and completes the proof of the theorem.
\end{proof}

We end this section with the following easy but useful lemma.

\bl\label{lemma:power-sum}
 Let $k\geq 2$, $\alpha>1$ and $0<\lambda, \delta\leq 1$. If $L>0$ and $0\leq t_i<\delta L$ are such that
 \begin{equation*}
  \lambda\sum_{i=1}^\infty t_i\leq L
 \end{equation*}
 then
 \begin{equation*}
  \sum_{i=1}^\infty I_{k,\alpha}(t_i) 
   \leq \frac{2(1+\delta\lambda)}{\lambda}\max\left\{(2\delta)^{\frac{1}{k-1}},(2\delta)^{\frac{1}{\alpha-1}}\right\}I_{k,\alpha}(L).
 \end{equation*}
\el
\begin{proof}
 Pick finitely many integer numbers $0=:m_0<m_1<m_2<\dots<m_{j_0}$ with the property that
 \begin{equation*}
  \delta L< t_{m_{i-1}+1}+\dots+t_{m_i}<2\delta L
 \end{equation*}
 for each $i=1,\dots,j_0$ and
 \begin{equation*}
  \sum_{n=m_{j_0}+1}^\infty t_n\leq \delta L.
 \end{equation*}
 Then $j_0\leq \frac{1}{\lambda\delta}$ and hence
 \begin{equation*}
  \begin{split}
   \sum_{i=1}^\infty I_{k,\alpha}(t_i) &\leq \sum_{i=1}^{j_0}I_{k,\alpha}(t_{m_{i-1}+1}+\dots+t_{m_i}) + I_{k,\alpha}\left(\sum_{n=m_{j_0}+1}^\infty t_n\right)\\
        &\leq \frac{1}{\lambda\delta}I_{k,\alpha}(2\delta L) + I_{k,\alpha}(\delta L)\\
        &\leq \frac{2(1+\delta\lambda)}{\lambda}\max\left\{(2\delta)^{\frac{1}{k-1}},(2\delta)^{\frac{1}{\alpha-1}}\right\} I_{k,\alpha}(L).
  \end{split}
 \end{equation*}
 We note that the inequality in the first line above is a consequence of (vii) of \lemref{lemma:trivial-properties-functions}.
\end{proof}

\section{Proof of the compactness theorem}

Let $X, X', X''$ be complete metric spaces with the property that $X$ isometrically embeds into $X'$ and $X'$ into $X''$. Let $k\geq 1$ and suppose $X''$ admits an isoperimetric inequality of Euclidean type for $\intcurr_k(X'')$ with some constant $D_k$. If $k\geq 2$ then suppose furthermore that $X'$ admits an isoperimetric inequality of Euclidean type for $\intcurr_{k-1}(X')$ with some constant $D_{k-1}$.

Choose a sequence $\frac{1}{2}>\lambda_1>\lambda_2>\dots >0$ of numbers satisfying
\begin{equation*}
 \prod_{i=1}^\infty\frac{1-\lambda_i}{1+\lambda_i} \geq \frac{1}{2}.
\end{equation*}
If $k=1$ set $\gamma_i:= \frac{1}{2}$ for all $i\geq 1$, if $k\geq 2$, let $\gamma_i= \gamma(\lambda_i)$ be the constant of \thmref{theorem:suitable-decomposition} corresponding to $\lambda_i$.

\bl\label{lemma:special-decomp-compactness}
Given a sequence $\frac{1}{2}>\delta_1>\delta_2>\dots>0$ and $T\in\intcurr_k(X)$ there exist sequences $T^i, S^i\in\intcurr_k(X')$ such that for all $m\in\N$
\begin{equation*}
 T= T^1+\dots+T^m+S^m
\end{equation*}
and the following properties hold:
\begin{enumerate}
 \item $\bdry T^1=\bdry T$ and $\bdry T^{j+1}=\bdry S^j=0$ for all $j\geq 1$;
 \item $\mass{T^1} +\dots+\mass{T^m}+\mass{S^m} \leq 2\mass{T}$;
 \item For every $1\leq j \leq m-1$
  \begin{equation*}
    \fillvol_{X''}(T^{j+1}+\dots+T^m+S^m)\leq 108D_k\delta_j\mass{T}^{\frac{k+1}{k}};
  \end{equation*}
 \item For every $j\geq 1$ and all $x\in\spt T^j$
  \begin{equation}\label{equation:growth-cptness-thm}
   \|T^j\|(B(x,r))\geq \frac{\gamma_j}{2\cdot25^k}\cdot r^k
  \end{equation}
  whenever  $0\leq r \leq \min\{5\delta_j\mass{T^j}^{\frac{1}{k}}, \dist(x,\spt\bdry T^j)\}$.
\end{enumerate}
\el

Here we use the convention that $\dist(x,\emptyset)=\infty$. In particular,  so \eqref{equation:growth-cptness-thm} holds for all $0\leq r \leq 5\delta_j\mass{T^j}^{\frac{1}{k}}$ whenever $j\geq 2$. The proof will show that if $k=1$ we may replace (ii) by
\begin{equation*}
 \mass{T^1} +\dots+\mass{T^m}+\mass{S^m} = \mass{T}.
\end{equation*}

\begin{proof}
 Fix $\delta, \lambda\in (0,1/2)$ and let $T'\in\intcurr_k(X')$. Suppose $T'= R+\sum_{i=1}^\infty T_i$ is a decomposition as in 
 \thmref{theorem:suitable-decomposition}. Setting $T^1:= R$ and $S^1:= \sum_{i=1}^\infty T_i$ we obtain
 \begin{equation*}
  \mass{T^1} + \frac{1-\lambda}{1+\lambda}\mass{S^1}\leq \mass{T'},
 \end{equation*}
 or $\mass{T^1} + \mass{S^1} = \mass{T'}$ if $k=1$.
By \lemref{lemma:power-sum}, we obtain
 \begin{equation*}
 \begin{split}
  \fillvol_{X''}(S^1)
 \leq \frac{(2+2\lambda)^{\frac{k+1}{k}}\gamma^{\frac{1}{k}}[1+(1-\lambda)\delta^k\gamma]}{1-\lambda}\delta D_k\mass{T'}^{\frac{k+1}{k}}
   \leq 27D_k\delta \mass{T'}^{\frac{k+1}{k}}.
 \end{split}
 \end{equation*}
 Using such a decomposition procedure successively with $\delta_i, \lambda_i$ we obtain a decomposition with the desired properties.
\end{proof}

\bp\label{proposition:extended-gromov-compactness-sets}
Let  $(X_n, d_n)$ be a sequence of metric spaces and, for each $n\in\N$, subsets
\begin{equation*}
 B_n^1\subset B_n^2\subset B_n^3\subset\dots\subset X_n.
\end{equation*}
If for every $i\in\N$ the sequence $(B_n^i, d_n)$ is uniformly compact then there exist a complete metric space $Z$, a subsequence $X_{n_m}$,
isometric embeddings $\varphi_m: X_{n_m}\hookrightarrow Z$ and compact subsets $Y^1\subset Y^2\subset\dots\subset Z$ and \begin{equation*}
\varphi_m(B_{n_m}^i)\subset Y^i\quad\text{ for all $m\in\N$ and $i\in\N$.}
\end{equation*}
\ep

\begin{proof}
 Choose for all $n,i\in\N$ a maximally $2^{-i}$-separated subset
 \begin{equation*}
 \{x(n,i,1), x(n,i,2),\dots, x(n,i,m(n,i))\}\subset B_n^i.
 \end{equation*}
 Clearly, $\sup_{n}m(n,i)<\infty$ for fixed $i\in\N$. After passage to a diagonal subsequence we may assume without loss of generality that for each $i_0\in\N$ we have $m(n,i_0) = m_{i_0}$ for every $n\geq i_0$ and that
 \begin{equation}\label{equation:distance-compare}
  |d_n(x(n,i,k), x(n, i', k')) - d_{n+1}(x(n+1, i, k), x(n+1, i', k'))|\leq \frac{1}{2^{n+1}}
 \end{equation} 
 for all $n\geq i_0$, $1\leq i,i'\leq i_0$, $1\leq k\leq m_i$ and $1\leq k'\leq m_{i'}$.
 We define a metric $d^n$ on the disjoint union $X^n:= X_1\sqcup X_2\sqcup\dots\sqcup X_n$ iteratively as follows. For $n=1$ simply set $d^1:= d_1$. If $n\geq 2$ and if the metric $d^{n-1}$ on $X^{n-1}$ has already been defined then we let $d^n$ be the unique metric on $X^n=X^{n-1}\sqcup X_n$ which equals $d^{n-1}$ on $X^{n-1}$, equals $d_n$ on $X_n$ and which, for $x\in X_n$ and $x'\in X^{n-1}$, is given by
 \begin{equation*}
  d^n(x,x'):= \frac{1}{2^n} + \min\{d_n(x,x(n,i,k)) + d^{n-1}(x', x(n-1, i, k))\}
 \end{equation*}
 where the minimum is taken over all $i,k$ satisfying $1\leq i\leq n-1$ and $1\leq k\leq m_i$. Using \eqref{equation:distance-compare} it is trivial to check that $d^n$ is indeed a metric. 
Finally, define $Z$ to be the completion of $\cup X^n$ with the metric
 coming from $d^n$. Clearly, the natural inclusion from $X_n$ to $Z$ is isometric for every $n$. We now show that for fixed $i$ the set $B^i:= \bigcup_{n=1}^\infty B_n^i$, as a subset of $Z$,
is pre-compact.
For this, note first that
 for each $i$, all  $n>n'\geq i$ and $1\leq k\leq m_i$
 \begin{equation*}
  d^n(x(n,i,k), x(n',i,k))\leq \sum_{l=n'}^{n-1} d^{l+1}(x(l,i,k),x(l+1,i,k))\leq \sum_{l=n'}^{n-1}2^{-(l+1)}\leq \frac{1}{2^{n'}}.
 \end{equation*}
 It follows that $B_n^i$ lies in the $(2^{1-i}+2^{-n'})$-neighborhood of $B_{n'}^i$. 
 Let $\varepsilon>0$ and choose $n'\geq i$ so large that $2^{-n'}\leq \varepsilon/8$. Then $\bigcup_{n\geq n'+1}B_n^{n'}$ is contained in the
 $\frac{\varepsilon}{2}$-neighborhood of $B_{n'+1}^{n'}$. Since $B_n^i\subset B_n^{n'}$ for all $n$, we can cover $B^i$ by finitely many balls of radius $\varepsilon$. This proves pre-compactness of $B^i$. Thus, $Y^i:=\bar{B^i}$ is compact and this concludes the proof.
\end{proof}

We are ready for the proof of the compactness theorem.

\begin{proof}[Proof of \thmref{theorem:general-compactness-fillvol}]
 Let $D_{k-2}$, $D_{k-1}$ and $D_k$ be suitable constants and $X_n\subset X'_n\subset X''_n$ complete metric spaces with the following properties: Firstly, $X''_n$ admits an isoperimetric inequality of Euclidean type for $\intcurr_k(X''_n)$ with constant $D_k$. Secondly, if $k\geq 2$, every closed ball $B\subset X'_n$ admits an isoperimetric inequality of Euclidean type for $\intcurr_{k-1}(B)$ with constant $D_{k-1}$ and, in case $k\geq 3$, also one for $\intcurr_{k-2}(B)$ with constant $D_{k-2}$. A possible choice would be $X''_n = X'_n = L^\infty(X_n)$.

Fix sequences $1=j_1<j_2<j_3<\dots$ of integers and $\frac{1}{2}>\delta_1>\delta_2>\dots>0$ of real numbers satisfying
\begin{equation*}
 \Delta:=\sum_{i=1}^\infty \delta_i<\infty.
\end{equation*}
Let $T_n\in\intcurr_k(X_n)$ be as in the hypothesis of the theorem and choose closed balls $B_n\subset X'_n$ of radius $D$ such that $\spt T_n\subset B_n$. If $k=1$ then set $U_n^i:= 0\in\intcurr_1(X'_n)$. If $k\geq 2$ then let 
\begin{equation*}
 \bdry T_n =  W_n^1+\dots+W_n^m+V_n^m
\end{equation*}
be decompositions as in \lemref{lemma:special-decomp-compactness} with $U_n^i, V_n^m\in\intcurr_{k-1}(X'_n)$. We may assume without loss of generality that for each $n\in\N$, all $V_n^m$ and $W_n^i$ are supported in $B_n$.
Set $\tilde{W}_n^{j_n+1}:= V_n^{j_n} = W_n^{j_n+1} + V_n^{j_n+1}$ and, for $1\leq i\leq j_n$, set $\tilde{W}_n^i:= W_n^i$. By \propref{proposition:isop-ineq-growth-estimate} there exists $U_n^i\in\intcurr_k(B_n)$ with $\bdry U_n^i = \tilde{W}_n^i$ and $\mass{U_n^i}\leq 2\fillvol(\tilde{W}_n^i)$
and
\begin{equation}\label{equation:mass-growth-bdry-decomp}
 \|U_n^i\|(B(x,r))\geq \frac{r^k}{(3D_{k-1})^{k-1}k^k}
\end{equation}
for all $x\in\spt U_n^i$ and every $0\leq r \leq \dist(x,\spt\tilde{W}_n^i)$.
In particular, we have
\begin{equation*}
 \mass{U_n^1}\leq 2\fillvol(W_n^1)\leq 2\cdot 2^{\frac{k}{k-1}}D_{k-1}\mass{\bdry T_n}^{\frac{k}{k-1}}
   \leq 2\cdot 2^{\frac{k}{k-1}}D_{k-1}C^{\frac{k}{k-1}}
\end{equation*}
and, for $2\leq i\leq j_n$,
\begin{equation*}
 \mass{U_n^i} 
\leq 2\left(\fillvol(W_n^i+W_n^{i+1}+V_n^{i+1})+\fillvol(W_n^{i+1}+V_n^{i+1})\right)
\leq 216D_{k-1}(\delta_{i-1}+\delta_i)C^{\frac{k}{k-1}}
\end{equation*}
and finally
\begin{equation*}
 \mass{U_n^{j_n+1}} \leq 2\fillvol(W_n^{j_n+1}+V_n^{j_n+1})
   \leq 216D_{k-1}\delta_{j_n}C^{\frac{k}{k-1}}.
\end{equation*}
From this it follows that
\begin{equation}\label{equation:mass-Ui}
 \sum_{i=1}^{j_n+1}\mass{U_n^i}\leq 2D_{k-1}C^{\frac{k}{k-1}}\left(2^{\frac{k}{k-1}}+216\Delta\right)
\end{equation}
and, for $1\leq L\leq j_n$,
\begin{equation}\label{equation:mass-Ui-from-L}
 \sum_{i=L+1}^{j_n+1}\mass{U_n^i}\leq 432 D_{k-1}C^{\frac{k}{k-1}}\sum_{i=L}^\infty\delta_i.
\end{equation}
Note that the right hand side of \eqref{equation:mass-Ui-from-L} tends to $0$ as $L\to\infty$.
Set 
\begin{equation*}
 \tilde{T}_n:= T_n - \sum_{i=1}^{j_n+1} U_n^i
\end{equation*}
and note that $\tilde{T_n}\in\intcurr_k(B_n)$. If $k=1$ then $\bdry\tilde{T}_n = T_n$; if $k\geq 2$, then $\bdry\tilde{T}_n = 0$ and, by \eqref{equation:mass-Ui},
\begin{equation*}
 \mass{\tilde{T}_n}\leq C + 2D_{k-1}C^{\frac{k}{k-1}}\left(2^{\frac{k}{k-1}}+216\Delta\right).
\end{equation*}
Let $\tilde{T}_n =  T_n^1+\dots+T_n^{j_n}+S_n^{j_n}$
be a decomposition as in \lemref{lemma:special-decomp-compactness} with $T_n^i, S_n^{j_n}\in\intcurr_k(X'_n)$. We may again assume without loss of generality that for each $n\in\N$, $S_n^{j_n}$ 
and all $T_n^m$ are supported in $B_n$. For each $i\in\N$ define
\begin{equation*}
 B_n^i:= \bigcup_{\nu=1}^{\min\{i, j_n\}}\left(\spt T_n^\nu \cup \spt U_n^\nu\right).
\end{equation*}
It follows that, for fixed $i$, the sequence $(B_n^i)$ is uniformly compact. This follows from the volume growth properties of $\tilde{W}_n^\nu$, $U_n^\nu$ and $T_n^\nu$ for $\nu=1,\dots,\min\{i, j_n\}$, see \eqref{equation:growth-cptness-thm} and \eqref{equation:mass-growth-bdry-decomp}, and from the fact that $\bdry\tilde{T}_n=0$ if $k\geq 2$, or that $\spt( \bdry\tilde{T}_n)$ consists of a uniformly bounded number of points if $k=1$.
Hence, by \propref{proposition:extended-gromov-compactness-sets}, 
there exist a complete metric space $Z$, isometric embeddings $\varphi_m: X''_{n_m}\hookrightarrow Z$ of a subsequence $X''_{n_m}$ and compact subsets $Y^i\subset Z$ with $Y^1\subset Y^2\subset\dots$ and 
\begin{equation*}
\varphi_m(B_{n_m}^i)\subset Y^i\quad\text{ for all $m\in\N$ and $i\in\N$.}
\end{equation*}
After possibly replacing $Z$ by $L^\infty(Z)$ we may assume without loss of generality that $Z$ is quasiconvex and admits a cone type inequality for $\intcurr_l(Z)$, for $l=1,\dots, k$, in the sense of \cite{Wenger-flatconv}.
After possibly choosing a diagonal sequence we may assume by the compactness and closure theorems for integral currents \cite{Ambr-Kirch-curr} that for each $i$ fixed, 
$\varphi_{m\#}T_{n_m}^i\rightharpoonup T^i$ and $\varphi_{m\#}U_{n_m}^i\rightharpoonup U^i$ as $m\to\infty$ for some $T^i, U^i\in\intcurr_k(Z)$.
For every $L\geq 1$ we have, by \eqref{equation:mass-Ui},
\begin{equation}\label{equation:limit-Ui-upto-L}
 \sum_{i=1}^L\mass{U^i} \leq \liminf_{m\to\infty}\sum_{i=1}^L\mass{U_{n_m}^i}\leq 2D_{k-1}C^{\frac{k}{k-1}}\left(2^{\frac{k}{k-1}}+ 216\Delta\right)
\end{equation}
and analogously
\begin{equation}\label{equation:limit-Ti-upto-L}
 \sum_{i=1}^L\mass{T^i}\leq 2C+4D_{k-1}C^{\frac{k}{k-1}}\left(2^{\frac{k}{k-1}}+216\Delta\right),
\end{equation}
furthermore
\begin{equation*}
 \sum_{i=1}^L\mass{\bdry U^i} \leq 2 C.
\end{equation*}
It follows that $T:= \sum_{i=1}^\infty (T^i+U^i)\in\intcurr_k(Z)$.
We finally show that
\begin{equation*}
 \flatnorm(T-\varphi_{m\#}T_{n_m}) \to 0\quad\text{as $m\to\infty$.}
\end{equation*}
For a fixed $L\in\N$ and large enough $m$ we write
\begin{equation*}
 \begin{split}
  T-\varphi_{m\#}T_{n_m} &=  \left[T-\sum_{i=1}^L(T^i+U^i)\right] + \sum_{i=1}^L(T^i- \varphi_{m\#}T_{n_m}^i) + \sum_{i=1}^L(U^i- \varphi_{m\#}U_{n_m}^i)\\
&\;\;\; - \varphi_{m\#}(T_{n_m}^{L+1}+\dots+T_{n_m}^{j_{n_m}} + S_n^{j_{n_m}})
   - \left[\sum_{i=L}^{j_{n_m}}\varphi_{m\#}U_{n_m}^{i+1}\right].
 \end{split}
\end{equation*}
By \eqref{equation:mass-Ui-from-L}, \eqref{equation:limit-Ui-upto-L} and \eqref{equation:limit-Ti-upto-L}, the expressions in the first and last brackets converge to $0$ in mass as $L\to\infty$ and therefore also in the flat norm. Each term in the two sums in the middle converges weakly to $0$ as $m\to\infty$. Since $Z$ is quasiconvex and admits cone type inequalities it follows from \cite{Wenger-flatconv} that their flat norms also converge to $0$. Finally, the filling volume of the remaining expression is arbitrary small when $L$ is large, by (iii) of \lemref{lemma:special-decomp-compactness}. Letting first $m$ and then $L$ tend to infinity we conclude the proof. As for the last statement in the theorem, it is enough to note that, by \cite{Wenger-flatconv}, weak convergence, flat convergence and convergence with respect to filling volume are equivalent in $Z$.
\end{proof}

\section{Uniqueness}\label{Section:Uniqueness}

\thmref{theorem:uniqueness-intro} is a consequence of the following result .

\bt\label{theorem:uniqueness-last-section}
 Let $Z, Z'$ be complete metric spaces, $k\geq 0$, and $T\in\intcurr_k(Z)$, $T'\in\intcurr_k(Z')$ integral currents whose supports have finite diameter. Suppose there exist  complete metric spaces $Z_n$, isometric embeddings $\varrho_n: \spt T\hookrightarrow Z_n$ and $\varrho'_n:\spt T'\hookrightarrow Z_n$ such that
 \begin{equation*}
  \flatnorm(\varrho_{n\#}T - \varrho'_{n\#}T')\to 0.
 \end{equation*}
 If $T\not=0$ then $T'\not=0$ and there exists an isometry $\Psi: \spt T \to \spt T'$ with $\Psi_\#T=T'$.
\et

\begin{proof}
 Suppose $T\not=0$.
 By Lemma 2.9 of \cite{Ambr-Kirch-curr} there exist compact sets $C_1\subset C_2\subset\dots\subset Z$ and $C'_1\subset C'_2\subset\dots\subset Z'$ so that $\|T\|(Z\backslash C_i)\leq 2^{-i}$ and $\|T'\|(Z'\backslash C'_i)\leq 2^{-i}$ for all $i\geq 1$. For $i, n\in\N$ set
\begin{equation*}
 B_n^i:=\varrho_n(C_i)\cup\varrho'_n(C'_i)
\end{equation*}
and note that $B_n^1\subset B_n^2\subset\dots\subset Z_n$. We claim that for fixed $i\in\N$, the sequence $(B_n^i)$ has uniformly bounded diameter. 
Suppose, in the contrary, that there exist $R>0$ and $z_n, \bar{z}_n\in Z_n$ satisfying $d_{Z_n}(z_n, \bar{z}_n)\to\infty$ and 
\begin{equation*}
\varrho_n(\spt T)\subset B(z_n, R)\quad\text{and}\quad \varrho'_n(\spt T')\subset B(\bar{z}_n, R).
\end{equation*}

By \propref{proposition:isop-ineq-growth-estimate} there exist $U_n\in\intcurr_k(Z_n)$ and $V_n\in\intcurr_{k+1}(Z_n)$ with 
\begin{equation*}
 \varrho_{n\#}T - \varrho'_{n\#}T' = U_n + \bdry V_n,
\end{equation*}
and such that $\mass{U_n}, \mass{V_n} \to 0$. After possibly replacing $Z_n$ by $l^\infty(Z_n)$ we may assume that $U_n$ and $V_n$ have the volume growth property of \propref{proposition:isop-ineq-growth-estimate}. In particular, it follows that for $R'>0$ large enough, $U_n$ and $V_n$ have support in $B(z_n, R'/2)\cup B(\bar{z}_n, R'/2)$. The slicing theorem in \cite{Ambr-Kirch-curr} then implies that 
\begin{equation*}
 \varrho_{\#n}T = (\varrho_{n\#}T - \varrho'_{n\#}T')\rstr B(z_n, R') = U'_n + \bdry V'_n,
\end{equation*}
for all sufficiently large $n$. Here we abbreviated $U'_n:=  U_n\rstr B(z_n, R')$ and $V'_n:= V_n\rstr B(z_n, R')$. Now, let $\iota: Z\hookrightarrow l^\infty(Z)$ be a Kuratowski embedding. For each $n$ choose a $1$-Lipschitz extension $\eta_n: Z_n\to l^\infty(Z)$ of $\iota\circ\varrho_n^{-1}$. Then 
\begin{equation*}
 \iota_\#T = \eta_{n\#}U'_n + \bdry(\eta_{n\#}V'_n).
\end{equation*}
Since $\mass{U'_n}$ and $\mass{V'_n}$ tend to $0$ it follows that $\flatnorm(T)= 0$ in $l^\infty(Z)$ and thus $T=0$. Since we assumed that $T\not=0$ this proves our claim on $(B_n^i)$.  It thus follows that for fixed $i\in\N$, the sequence $(B_n^i)$ is uniformly compact. After passing to a suitable subsequence we may assume by \propref{proposition:extended-gromov-compactness-sets} that there exist a complete metric space $Z''$, isometric embeddings $\psi_n: Z_n\hookrightarrow Z''$ and compact subsets $Y^1\subset Y^2\subset\dots\subset Z''$ such that
\begin{equation*}
 \psi_n(B_n^i)\subset Y^i\quad\text{for all $n\in\N$ and $i\in\N$.}
\end{equation*}
Consider the isometric embedding $\psi_n\circ\varrho_n: \spt T\hookrightarrow Z''$ and first note that $\psi_n\circ\varrho_n(C_i)\subset Y^i$. Therefore, after passing to a subsequence, $\psi_n\circ\varrho_n$ converges pointwise to an isometric embedding $\psi: \spt T\hookrightarrow Z''$, uniformly on each $C_i$. We now claim that $(\psi_n\circ\varrho_n)_\#T$ converges weakly to $\psi_\#T$. Indeed, for $f, \pi_j\in\lip(Z'')$ with $f$ bounded, we have
\begin{equation*}
 \begin{split}
 |(\psi_n&\circ\varrho_n)_\#T(f, \pi) - \psi_\#T(f,\pi)|\\
 &\leq |T(f\circ\psi_n\circ\varrho_n - f\circ\psi, \pi\circ\psi_n\circ\varrho_n)| 
     + |T(f\circ\psi, \pi\circ\psi_n\circ\varrho_n) - T(f\circ\psi, \pi\circ\psi)|.
  \end{split}
\end{equation*}
The first term converges to $0$ as $n\to\infty$ by the choice of $C_i$ and the definition of mass. The second term converges to $0$ by the continuity property of currents. This shows that $(\psi_n\circ\varrho_n)_\#T$ converges weakly to $\psi_\#T$, as claimed. Analogously, after passing to a subsequence, we may assume that $\psi_n\circ\varrho'_n$ converges pointwise to an isometric embedding $\psi': \spt T'\hookrightarrow Z''$, uniformly on each $C'_i$ and that $(\psi_n\circ\varrho'_n)_\#T'$ converges weakly to $\psi'_\#T'$. Finally, we have $\psi_\#T = \psi'_\#T'$ because in the following equality
\begin{equation*}
 \psi_\#T - \psi'_\#T' = \left[\psi_\#T - (\psi_n\circ\varrho_n)_\#T\right] + \left[(\psi_n\circ\varrho'_n)_\#T' - \psi'_\#T'\right] + \psi_{n\#}(\varrho_{n\#}T - \varphi'_{n\#}T')
\end{equation*}
all terms converge weakly to $0$. Indeed, for the first and second term, this was shown above. For the third term this follows because
\begin{equation*}
 \flatnorm(\psi_{n\#}(\varrho_{n\#}T - \varrho'_{n\#}T'))\leq \flatnorm(\varrho_{n\#}T - \varrho'_{n\#}T') \to 0.
\end{equation*}
This shows that $\psi_\# T = \psi'_\#T'$ and hence 
\begin{equation*}
 \psi(\spt T) = \spt(\psi_\#T) = \psi'(\spt T').
 \end{equation*}
 It follows that $\Psi:= \psi'^{-1}\circ\psi: \spt T\to\spt T'$ is an isometry which satisfies
$\Psi_\#T=T'$.
\end{proof}

We are ready for the proof of the uniqueness result stated in the introduction.

\begin{proof}[{Proof of \thmref{theorem:uniqueness-intro}}]
 Let $(Z, d)$ and $(Z', d')$ be metric spaces as in the hypothesis. Let $d_n$ be the unique pseudo-metric on $Y:= Z\sqcup Z'$ agreeing with $d$ on $Z$, with $d'$ on $Z'$ and so that
 \begin{equation*}
  d_n(z,z'):= \inf\{d(z,\varphi_n(x)) + d'(\varphi'_n(x), z'): x\in\spt T_n\}
 \end{equation*}
 for all $z\in Z$ and $z'\in Z'$.
 Let $(Z_n, d_n)$ be the completion of the metric space associated with $(Y, d_n)$ and let $\varrho_n: \spt T\hookrightarrow Z_n$ and $\varrho'_n: \spt T' \hookrightarrow Z_n$ be the natural inclusions. Clearly, $\varrho_n$ and $\varrho'_n$ are isometric embeddings and satisfy
 \begin{equation*}
  \flatnorm(\varrho_{n\#}T - \varrho'_{n\#}T')\to 0.
 \end{equation*}
 The theorem now follows from \thmref{theorem:uniqueness-last-section}.
\end{proof}

\end{document}